\documentclass[12pt]{amsart}
\usepackage{amsmath, amssymb}

\newtheorem{theorem}{\sc Theorem}[section]
\newtheorem{lemma}[theorem]{\sc Lemma}
\newtheorem{proposition}[theorem]{\sc Proposition}

\newtheorem{hypothesis}[theorem]{\sc Hypothesis}

\newtheorem{example}[theorem]{\sc Example}

\newtheorem*{thmMT}{\sc Theorem \ref{MT}}
\newtheorem*{thmEbfg}{\sc Theorem \ref{Ebfg}}

\makeatletter
\@namedef{subjclassname@2020}{\textup{2020} Mathematics Subject Classification}
\makeatother

\usepackage{abstract}

\title[Finite conjugacy classes]{Profinite groups with finite conjugacy classes of commutators}

\author[E. Detomi]{Eloisa Detomi}
\address{Eloisa Detomi: Dipartimento di Matematica \lq\lq Tullio Levi-Civita\rq\rq, Universit\`a degli Studi di Padova, Via Trieste 63, 35121 Padova, Italy} 
\email{eloisa.detomi@unipd.it}

\author[D. Senise]{D\'ebora Senise} 
\address{D\'ebora Senise: Department of Mathematics, University of Brasilia, Brasilia, DF, Brazil}
\email{D.F.P.Senise@mat.unb.br}

\author[P. Shumyatsky]{Pavel Shumyatsky} 
\address{Pavel Shumyatsky: Department of Mathematics, University of Brasilia, Brasilia, DF, Brazil}
\email{pavel@unb.br}


\begin{document}

\maketitle

\begin{abstract} It was shown in a recent work by Senise and Shumyatsky that if $G$ is a finite group such that $|x^G|\leq m$ whenever $x$ is a commutator of prime power order, then $G''$ has $m$-bounded order. In this paper, we examine the following question:  given a profinite group $G$ such that $|x^G|$ is finite whenever $x$ is a commutator of prime power order, does it follow that $G''$ is finite? We show that, in general, the answer to this is negative. The main result obtained in this paper is the following theorem.

\noindent Let $G$ be a profinite group in which $|x^G|$ is finite whenever $x$ is a commutator of prime power order. Assume that $\pi(G)$ is finite. Then
    \begin{enumerate}
\item The second commutator subgroup $G''$ is finite.
\item $G$ has an open normal subgroup of nilpotency class at most $3$.
    \end{enumerate}
\end{abstract}


\section{Introduction}

Given a group $G$ and an element $x\in G$, we write $x^G$ for the conjugacy class of $x$. 
In 1954, B. H. Neumann proved that if all conjugacy classes in a group $G$ are finite and have size at most $m$, then the order of the commutator subgroup $G'$ is finite and $m$-bounded \cite{neumann}. We say that a quantity is $(a,b,c\dots)$-bounded if it is bounded by a constant depending only on the parameters $a,b,c\dots$. 
 There are several recent publications in which Neumann's theorem has been extended in various directions (see in particular \cite{as,dms,dieshu}). The following theorem collects some relevant results on groups $G$ in which $|x^G| \le m$ whenever $x$ is a commutator, that is, $x=a^{-1}b^{-1}ab$ for some $a,b\in G$: 
Part 1 of the theorem was established in \cite{dieshu}, while Part 2 is immediate from the main result of \cite{ebeshu}.

\begin{theorem} \label{known1}
Let $m$ be a positive integer and $G$ a group in which $|x^G|\leq m$ whenever $x$ is a commutator. Then 

\begin{enumerate}
\item The second commutator subgroup $G''$ is finite of $m$-bounded order.
\item $G$ has a normal subgroup $H$ such that the index $[G:H]$ and the order of $\gamma_4(H)$ are both $m$-bounded.
\end{enumerate}
\end{theorem}

In the case of finite groups, Neumann's theorem admits the following generalizations, as shown in \cite{ds25}.

\begin{theorem} \label{ds}
     Let $m$ be a positive integer and $G$ a finite group. Then 

    \begin{enumerate}
        \item If $|x^G|\leq m$ whenever $x$ is an element of prime power order, then $G'$ has $m$-bounded order.
        \item If $|x^G|\leq m$ whenever $x$ is a commutator of prime power order, then $G''$ has $m$-bounded order.
    \end{enumerate} 
\end{theorem}

In \cite{ds25}, the authors observed that Part 2 of Theorem \ref{known1} cannot be extended to commutators of prime power order in a similar manner. However, an analogous result holds if we also take into account the number of primes dividing the order of $G$.

{
\renewcommand{\thetheorem}{A}

\begin{theorem} \label{Ebfg}
 Let $m,s$ be positive integers and $G$ a finite group such that $|\pi(G)|=s$ and $|x^G| \le m$ whenever $x$ is a commutator of prime power order. Then $G$ has a normal subgroup $H$ such that the index $[G : H]$ and the order of $\gamma_4(H)$ are both $(s,m)$-bounded.
\end{theorem}
\addtocounter{theorem}{-1}
}

\renewcommand{\thetheorem}{\arabic{section}.\arabic{theorem}}

Throughout, $\pi(G)$ denotes the set of primes $p$ for which $G$ has a nontrivial Sylow $p$-subgroup. In the finite case, $\pi(G)$ is precisely the set of primes dividing $|G|$.

We remark that the last two results can be extended to profinite groups by an inverse limit argument. 
Here we are interested in profinite groups $G$ in which $|x^G|$ is finite whenever $x$ is a commutator. The following results are known (see \cite{dms20} and \cite{dms23}).

\begin{theorem} \label{known}
    Let $G$ be a profinite group in which $|x^G|$ is finite whenever $x$ is a commutator. Then
    \begin{enumerate}
\item The second commutator subgroup $G''$ is finite.
\item $G$ has an open normal subgroup of nilpotency class at most $3$.
    \end{enumerate}
\end{theorem}

In this paper, we investigate whether Theorem~\ref{known} can be extended to profinite groups $G$ in which $|x^G|$ is finite whenever $x$ is a commutator of prime power order. Example~\ref{exam} shows that this is not possible in general. Although every such group is virtually prosoluble (see Proposition~\ref{Gprosol}), it need not be virtually nilpotent, and $G''$  need not be finite. However, as in Theorem~\ref{Ebfg}, the result does hold under the additional assumption that $\pi(G)$ is finite.

More precisely, we prove the following theorem.

\renewcommand{\thetheorem}{B}

\begin{theorem}  \label{MT}
 Let $G$ be a profinite group in which $|x^G|$ is finite whenever $x$ is a commutator of prime power order. Assume that $\pi(G)$ is finite. Then
    \begin{enumerate}
\item The second commutator subgroup $G''$ is finite.
\item $G$ has an open normal subgroup of nilpotency class at most $3$.
    \end{enumerate}
\end{theorem}

Recall that an element $x$ of a group $G$ is called a simple commutator of weight $n$ if there are $x_1,\dots,x_n\in G$ such that $x=[x_1,\dots,x_n]$. In view of Theorem \ref{MT}, it would be interesting to see if a similar result holds for profinite groups $G$ such that $\pi(G)$ is finite and $|x^G|$ is finite whenever $x$ is a simple commutator of weight $n$ of prime power order. So far, we have not been able to answer this question for $n\geq3$.

\renewcommand{\thetheorem}{\arabic{section}.\arabic{theorem}}


\section{Preliminaries}

Results on finite groups often admit a natural analog for profinite groups (see, for example, \cite{rz00} or \cite{wilson2}). Throughout the paper, we use certain profinite versions of facts on finite groups without explaining in detail how these can be deduced from the corresponding finite cases. Whenever this happens, the deduction follows by a standard inverse limit argument.  

Every homomorphism between profinite groups considered in this paper is continuous, and every subgroup of a profinite group is closed unless otherwise specified.

If $x$ is an element of a group $G$, we write $x^G$ for the conjugacy class of $x$ in $G$. We say that $x$ is an $FC$-element of $G$ if $|x^G|$ is finite. 
 We denote by $\Delta(G)$ the set of $FC$-elements of $G$, that is,
\[
\Delta(G) = \{ x \in G \mid |x^G| < \infty \}.
\]

The set $\Delta(G)$ is an abstract normal subgroup of $G$. Note that if $G$ is a profinite group, then $\Delta(G)$ need not be closed.

Recall that a pro-$p$ group is an inverse limit of finite $p$-groups, a pronilpotent group is an inverse limit of finite nilpotent groups, and a prosoluble group is an inverse limit of finite soluble groups. Profinite groups have Sylow $p$-subgroups and satisfy analogs of the Sylow theorems. 
\medskip

A difficulty that arises when working with commutators of prime power order is that the hypothesis may fail in a homomorphic image of the profinite group $G$. More precisely, if $Y(G)$ denotes the set of commutators of prime power order in $G$, and $N$ is a normal subgroup of $G$, then $Y(G/N)$ need not be the image of $Y(G)$ in $G/N$.

To circumvent this obstacle, we will work with the subset $X(G)$ of $Y(G)$, consisting of those commutators $[x,y]$ for which $x$ and $[x,y]$ belong to some Sylow subgroup $P$ of $G$, while $y$ is an arbitrary element of $G$. Throughout, we denote by $Syl(G)$ the set of Sylow subgroups of $G$. Thus,
 \[
 X(G)=\{[x,y]  \mid \text{there is } P \in Syl(G) \text{ such that }  x,[x,y]\in P, \  y\in G\}.
 \]

The main advantage of working with the set $X(G)$ rather than $Y(G)$ is that $X(G)$ is invariant under taking quotients. More precisely, the following lemma was proved in \cite[Lemma~3.1]{ds25} for finite groups, and the same proof applies to profinite groups.

\begin{lemma} \label{X(G)}
Let $G$ be a profinite group, and let $N$ be a normal subgroup of $G$. Then $ X(G)N/N=X(G/N)$. 
\end{lemma}

Let $G$ be a profinite group, $P$ a Sylow $p$-subgroup of $G$, and $N$ an open normal subgroup of $G$. Set $G_N = G/N$ and $P_N = PN/N$.  
It follows from the Focal Subgroup Theorem \cite[Theorem 7.3.4]{go} that
\[
\langle P_N \cap X(G_N) \rangle = P_N \cap G_N',
\]
hence $\langle X(G_N) \rangle = G_N'$. Consequently, $X(G)$ topologically generates $G'$. 
\medskip

We conclude this section with an example of a profinite group $G$ in which $|x^G|$ is finite whenever $x$ is a commutator of prime power order, while $G''$ is infinite. If $A$ and $B$ are groups, we denote by $B \wr A$ the wreath product of $B$ and $A$, that is, $B \wr A = (\prod_{|A|}B) \rtimes A$, where $A$ acts by permuting the factors of $(\prod_{|A|}B)$.

\begin{example} \label{exam}
    Let $\{ p_i \mid i \in \mathbb{N} \}$ be the set of all primes, and let $C_i$ be the cyclic group of order $p_i$. Set $P_i = (C_i \wr C_i) \wr C_i$. Since $P_i$ is finite, of course $|x^{P_i}|$ is finite whenever $x$ is a commutator of $P_i$. Moreover, $P_i'' \neq 1$ for every $i \in \mathbb{N}$. Consider the group $G = \prod_i P_i$.
    Every commutator of prime power order of $G$ must be a commutator of $P_i$ for some $i \in \mathbb{N}$. Thus, $G$ is a profinite group in which $|x^G|$ is finite whenever $x$ is a commutator of prime power order in $G$. On the other hand, $G''$ is infinite. 
\end{example}


\section{Proof of Theorem \ref{Ebfg}}

In this section, we prove Theorem \ref{Ebfg}, which we restate here for the reader's convenience. 

\begin{thmEbfg}
Let $m,s$ be positive integers and $G$ a finite group such that $|\pi(G)|=s$ and $|x^G| \le m$ whenever $x$ is a commutator of prime power order. Then $G$ has a normal subgroup $H$ such that the index $[G : H]$ and the order of $\gamma_4(H)$ are both $(s,m)$-bounded.
\end{thmEbfg}
  
\begin{proof}
    It follows from Theorem \ref{ds} (2) that $G''$ has $m$-bounded order. It suffices to prove the result for finite groups $G$ in which $|x^G| \le m$ whenever $x \in X(G)$.  By Lemma \ref{X(G)}, $X(G)$ is invariant under taking quotients. Thus, passing to the quotient $G/G''$, we may assume that $G$ is metabelian.  Let $F = F(G)$ be the Fitting subgroup of $G$. Write
\[
F = S_1 \times \cdots \times S_k,
\]
where $S_i$ is a Sylow $p_i$-subgroup of $F$ and $k\leq s$.

Let us show that $[G : F]$ is $(s,m)$-bounded. It suffices to prove that, for every Sylow $p$-subgroup $P$ of $G$, the quotient $P/O_p(G)$ has $(s,m)$-bounded order.

Suppose first that $O_p(G) = 1$, in which case $F$ is a $p'$-subgroup. Since $G$ is metabelian, $G' \le F$ and hence $P'= 1$. Thus $P$ is an abelian $p$-subgroup. For each $i = 1, \dots, k$, consider the subgroup $P S_i$. Since $P$ is abelian, $(P S_i)'\leq S_i$ and so every commutator of $PS_i$ has prime power order.

By applying Theorem \ref{known1} (2), there exists a normal subgroup $M_i$ of $P S_i$ such that $[P S_i: M_i]$ and $|\gamma_4(M_i)|$ are $m$-bounded. Let $C_i = C_{M_i}(\gamma_4(M_i))$. Note that $C_i$ is a normal subgroup of $P S_i$ of $m$-bounded index. 
Moreover, $\gamma_4(C_i)$ is central in $C_i$.  Thus $C_i$ is nilpotent, and so $C_i \le F(P S_i)$.

Hence, $F(P S_i)$ has $m$-bounded index in $P S_i$. Furthermore, $F(P S_i) = S_iC_P(S_i)$,
so the index $[P:C_P(S_i)]$ is $m$-bounded. 
Since this holds for each $i = 1, \dots, k$, 
the index $[P:C_P(F)]$ is $(s,m)$-bounded. As $G$ is soluble, $C_G(F) \le F$. Hence $C_P(F) = 1$ because $F$ is a $p'$-subgroup.  It follows that $P$ has $(s,m)$-bounded order. Removing the assumption that $O_p(G) = 1$, we conclude that $P/O_p(G)$ has $(s,m)$-bounded order. This shows that indeed $[G : F]$ is $(s,m)$-bounded.

Next, observe that every commutator of $S_i$ belongs to $X(G)$. 
Hence, by Theorem \ref{known1} (2), for each $i = 1, \dots, k$ there exists a normal subgroup $L_i$ of $S_i$ such that $[S_i:L_i]$ and $|\gamma_4(L_i)|$ are $m$-bounded. Set $L = \prod_{i=1}^{k} L_i$. 
We have that $[F:L]$ and $|\gamma_4(L)|$ are both $(s,m)$-bounded. Taking into account that $[G : F]$ is $(s,m)$-bounded, we remark that $[G : L]$ is $(s,m)$-bounded as well. Let $H$ be the normal core of $L$ in $G$. 
Then $H \trianglelefteq G$ and both $[G:H]$ and $|\gamma_4(H)|$ are $(s,m)$-bounded. 
The proof is complete.
\end{proof}

\section{$G$ is virtually nilpotent} 

The main goal of this section is to prove the following proposition.

\begin{proposition} \label{Gvirtnilp}
Let $G$ be a profinite group in which $\pi(G)$ is finite and $X(G) \subseteq \Delta(G)$. Then $G$ is virtually nilpotent.
\end{proposition}

Throughout this section, we work under the following hypothesis.

\begin{hypothesis} \label{Hyp}
Let $G$ be a profinite group in which $\pi(G)$ is finite and $X(G) \subseteq \Delta(G)$.
\end{hypothesis}

Of course, in every finite quotient $G/N$ of $G$, the elements of $X(G/N)$ are $FC$-elements.

\begin{lemma} \label{pronilp}
    Under Hypothesis \ref{Hyp}, if $H$ is a pronilpotent subgroup of $G$, then every commutator of $H$ belongs to $\Delta(G)$. 
\end{lemma}

\begin{proof}
        Let $k = |\pi(H)|$ and write 
\[
H = P_1 \times \dots \times P_k 
\]
    where $P_i$ is the Sylow $p_i$-subgroup of $H$. Every commutator of $H$ is of the form $[x,y]$ where $x = x_1 \cdots x_k$ 
    and $y = y_1 \cdots y_k$, with $x_i, y_i \in P_i$. Hence, $[x,y] = [x_1, y_1] \cdots [x_k, y_k]$. Since $[x_i,y_i] \in X(G) \subseteq \Delta(G) $ for $i = 1, \dots, k$, it follows that $[x,y]$ belongs to $\Delta(G)$ as well.   
\end{proof}

Recall that $\gamma_{\infty}(G)$ stands for the intersection of the terms of the lower central series of a group $G$. A profinite group $G$ is pronilpotent if and only if $\gamma_{\infty}(G)=1$. Each profinite group $G$ has a maximal pronilpotent normal subgroup, its Fitting subgroup $F(G)$. Set $F_0(G)=1$ and $F_{i+1}(G)/F_i(G)=F(G/F_i(G))$ for every integer $i\geq 0$. We say that the profinite group $G$ has (finite) Fitting height $h=h(G)$ if $G=F_h(G)$ and $h$ is the least integer with this property. Note that $G$ has finite Fitting height at most $h$ if and only if $G$ is an inverse limit of finite soluble groups of Fitting height at most $h$. A profinite group $G$ is metapronilpotent if and only if $h(G) \le 2$ or, equivalently, if and only if $\gamma_{\infty}(G)$ is pronilpotent.

Any prosoluble group $G$ has a Sylow basis (a family of pairwise  
permutable Sylow $p_i$-subgroups $P_i$ of $G$, exactly one for each prime
$p_i\in\pi(G)$), and any two Sylow bases are conjugate (see
\cite[Proposition~2.3.9]{rz00}). The basis normalizer (also known as the system
normalizer) of such a Sylow basis in $G$ is $T=\bigcap_i N_G(P_i)$. If $G$ is a
prosoluble group and $T$ is a basis normalizer in $G$, then $T$ is
pronilpotent and $G=\gamma_{\infty}(G)T$ (see \cite[Lemma~5.6]{r13}).

\begin{lemma} \label{[x,t]}
    Under Hypothesis \ref{Hyp}, if $G$ is prosoluble and $T$ is a basis normalizer in $G$, then for every $x \in G$ and $t \in T$ the commutator $[x, t]$ belongs to $\Delta(G)$.
\end{lemma}

\begin{proof}
Let $s = |\pi(G)|$ and let $\{ P_1, \dots , P_s \}$ be a Sylow basis of $G$ such that $T = \bigcap_{i=1}^s N_G(P_i)$.
Since $G = P_1P_2 \cdots P_s$, every element $x \in G$ can be written as a product $x = x_1x_2 \cdots x_s$ for some $x_i \in P_i$. Let $t \in T$ and write
\[ 
[x, t] = [x_1x_2 \cdots x_s, t] = [x_1, t]^{x_2 \cdots x_s}[x_2, t]^{x_3 \cdots x_s} \cdots [x_s, t].
\]

Since $T$ normalizes $P_i$, we have $[x_i, t] \in P_i$, and so $[x_i, t] \in X(G) \subseteq \Delta(G) $ for $i = 1, \dots, s$. As $\Delta(G)$ is normal in $G$, the commutator $[x,t]$ is a product of at most $s$ elements of $\Delta(G)$. Hence, $[x,t] \in \Delta(G)$ as well. 
\end{proof}

\begin{lemma} \label{metapronilp}
    Under Hypothesis \ref{Hyp}, if $G$ is a metapronilpotent group, then every commutator of $G$ belongs to $\Delta(G)$.
\end{lemma}

\begin{proof}
Since $G$ is metapronilpotent, and in particular prosoluble, we can write $G = NT$, where $T$ is a basis normalizer and $ N = \gamma_{\infty}(G).$ 
Moreover, in this case both $N$ and $T$ are pronilpotent. 
Let $x$ be a commutator of $G$, and write $x = [n_1t_1, n_2t_2] $
where $n_1, n_2 \in N$ and $t_1, t_2 \in T$. Then
\[
x = [n_1,t_2]^{t_1}[t_1,t_2][n_1,n_2]^{t_1t_2}[t_1,n_2]^{t_2}.
\]
By Lemma \ref{pronilp} and Lemma \ref{[x,t]}, all the commutators above belong to $\Delta(G)$, so $x \in \Delta(G)$ as well.
\end{proof}

The following lemmas do not require the assumption that $\pi(G)$ is finite.
\medskip

\begin{lemma} \label{simple}
 Let $G$ be a profinite group in which $X(G) \subseteq \Delta(G)$.
If $G$ is a Cartesian product of nonabelian finite simple groups,
then $G$ is finite.
\end{lemma}

\begin{proof} 
Let $G = \prod_{i \in I} S_i$, where each $S_i$ is a nonabelian finite simple group. 
Recall that simple groups have even order \cite{ft}. For every $i \in I$, choose a Sylow 2-subgroup $P_i$ of $S_i$ and elements $x_i\in P_i$ and $y_i\in S_i$ such that $1\neq [x_i,y_i]\in P_i$. 
 The existence of such elements follows from  the fact that 
\[
\langle X(S_i) \cap P_i \rangle = P_i \cap S_i'= P_i.
\]

Set $g=\prod_i [x_i,y_i]$. 
Then  $g=\left[\prod x_i,\prod y_i\right]$ and $\prod x_i$ belong to $\prod P_i$, which is a Sylow $2$-subgroup of $G$. Hence $g \in X(G) \subseteq \Delta(G)$. Thus $|g^G|$ is finite,  which forces $I$ to be finite. 
\end{proof}

Any finite group $H$ has a normal series such that each of its factors is either soluble or a direct product of nonabelian simple groups. The nonsoluble length of $H$,
denoted by $\lambda(H)$, is defined as the minimum number of nonsoluble factors in a
series of this kind. It is easy to see that the nonsoluble length $\lambda(H)$ is equal to the
least positive integer $l$ such that there is a series of characteristic subgroups
\[
1 = L_0 \leq R_0 < L_1 \leq R_1 < \cdots \leq R_l = H
\]
in which each quotient $L_i/R_{i-1}$ is a (nontrivial) direct product of nonabelian simple groups, and each quotient $R_i/L_i$ is soluble (possibly trivial).

It is natural to say that a profinite group $G$ has finite nonprosoluble length $l$ if $G$ has a normal series
\[
1 = L_0 \leq R_0 < L_1 \leq R_1 < \cdots \leq R_l = G
\]
in which each quotient $L_i/R_{i-1}$ is a (nontrivial) Cartesian product of nonabelian
finite simple groups, and each quotient $R_i/L_i$ is prosoluble (possibly trivial). In
particular, if for some positive integer $d$ all continuous finite quotients of a profinite group $G$ have nonsoluble length at most $d$, then $G$ has finite nonprosoluble length
at most $d$ (see e.g. \cite[Lemma 2]{wilson}).

\begin{proposition} \label{Gprosol}
 Let $G$ be a profinite group in which $X(G) \subseteq \Delta(G)$. Then $G$ is virtually prosoluble.
\end{proposition}

\begin{proof}
Suppose that $G$ is not prosoluble; thus, $G$ has a nontrivial Sylow $2$-subgroup $P$.
Since every commutator of $P$ is an element of $X(G) \subseteq \Delta(G)$, it follows from Theorem \ref{known} (2) that $P$ is virtually nilpotent. Hence,  $P$ is soluble, say of derived length $d$. 
It follows from \cite[Theorem 1.1]{ks15} that every finite image
of $G$ has nonsoluble length at most $d$, and hence $G$ also has nonprosoluble length at
most $d$.
 Let
\[
1 = L_0 \leq R_0 < L_1 \leq R_1 < \cdots \leq R_l = G
\]
be a normal series  of finite length $l \le d$ in which each section $L_i/R_{i-1}$ is a (nontrivial) Cartesian product of nonabelian finite simple groups, and each section $R_i/L_i$ is prosoluble
(possibly trivial). 
Let 
\[
H = \{ x \in G \mid [L_i,x] \leq R_{i-1} \text{ for all } i \} 
\]
 be the centralizer in $G$ of all the nonprosoluble sections of the series. By Lemma \ref{simple}, each 
 section $L_i/R_{i-1}$ is finite. 
 Hence, $[G : C_G(L_i/R_{i-1})]$ is finite for each $i = 1, \dots, l$. 
Since $l$ is finite as well, it follows that $H$ is open in $G$. 

Moreover, since the image of $H\cap L_i$ is central in $H/(R_{i-1}\cap H)$, it follows that $H$ is prosoluble. Hence, $G$ is virtually prosoluble, as required.
\end{proof}

\begin{lemma} \label{Fneq1}
 Let $G$ be a profinite group in which $X(G) \subseteq \Delta(G)$. If $G$ is a nontrivial prosoluble group, then $F(G) \neq 1$.
\end{lemma}

\begin{proof}
Assume by contradiction that $F(G) = 1$, and note that $X(G) \neq 1$. Indeed, since $ \langle X(G) \rangle = G'$, if $X(G) = 1$, then $G$ is abelian and $F(G) = G = 1$, which is not the case. Let $x$ be a nontrivial element of $X(G)$. As $X(G) \subseteq \Delta(G)$, the subgroup $K = \langle x^G \rangle$ is generated by finitely many conjugates of $x$, so
its centralizer has finite index in $G$. It follows that the center $Z(K)$ of $K$ has finite index in $K$. Moreover, $Z(K) \le F(G) = 1$. Therefore, $Z(K) = 1$, and consequently, $K$ is finite. As $G$ is prosoluble, $K$ is soluble and $F(K) \neq 1$. As $F(K) \le F(G)$, this
contradicts the hypothesis that $F(G)= 1$.
\end{proof}
\medskip

Now we are ready to prove Proposition \ref{Gvirtnilp}.
\medskip

\begin{proof}[Proof of Proposition \ref{Gvirtnilp}]
By Proposition \ref{Gprosol}, we may assume that $G$ is a nontrivial prosoluble group.
Hence, by Lemma \ref{Fneq1}, $F(G) \neq 1$.
It follows from Lemma \ref{metapronilp} and Theorem \ref{known} (2) that $F_2(G)$ is virtually nilpotent, that is, $F_2(G)/F(G)$ is finite. 
Let $G_F = G/F(G)$. Since $F(G_F)$ is finite, there exists an open normal subgroup $K_F$ of $G_F$ such that $K_F \cap F(G_F) = 1$.
Thus $F(K_F) = 1$. By Lemma \ref{Fneq1}, this forces $K_F = 1$. Therefore, $G_F$ is finite
and $G$ is virtually nilpotent.
\end{proof}
\medskip


\section{Proof of Theorem \ref{MT}} 

We will need the following lemma. 

\begin{lemma} \label{auxlemma}
    Let $G$ be a prosoluble group of finite Fitting height in which $\pi(G)$ is finite. Then $G$ has a series
\[
1 \le P_n \le P_nP_{n-1} \le \dots \le P_nP_{n-1}P_{n-2} \cdots P_1 = G,
\]
    where each $P_i$ is a pro-$p$ subgroup of $G$ for some $p \in \pi(G)$, and $P_i$ normalizes $P_j$ whenever $i \le j$.
\end{lemma}

\begin{proof}
    Set $s = |\pi(G)|$. Let $\{ P_1, \dots , P_s \}$ be a Sylow basis of $G$, and let $T = \bigcap_i N_G(P_i)$ 
    be its basis normalizer. Then $T$ is pronilpotent and $G = \gamma_{\infty}(G)T$. 
    Let $N_1 = \gamma_{\infty}(G)$ and $T_1 = T$. 
    Thus, 
    $\{ P_i \cap N_1 \}_i $ forms a Sylow basis of $N_1$, and 
    \[
    T_2 = \bigcap_{i=1}^s N_{N_1}(P_i \cap N_1)
    \]
    is the basis normalizer of $N_1$. Hence, $T_2$ is pronilpotent and $N_1 = \gamma_{\infty}(N_1)T_2$. 
    Since $N_1 \trianglelefteq G$, and $T_1$ normalizes each $P_i$, it follows that
    $T_1$ normalizes $T_2$. 
    Write $T_1 = Q_1 \times \dots \times Q_s$ 
    and $T_2 = R_1 \times \dots \times R_s$, where $Q_i$ and $R_i$ are Sylow $p_i$-subgroups of $T_1$ and $T_2$, respectively. The subgroup $T_2T_1$ therefore admits a series
    \[
    R_1 \le R_1R_2 \le \dots \le R_1 \cdots R_{s-1} \le T_2 \le T_2Q_1 \le \dots \le T_2T_1,
    \]
    satisfying the required conditions.
    
    Let $N_2 = \gamma_{\infty}(N_1)$. Since $N_1 = N_2T_2$, we have $G = N_2T_2T_1 $. 
    Similarly, set $N_2 = \gamma_{\infty}(N_2)T_3$, where $T_3$ is the basis normalizer of the Sylow basis $\{ P_i \cap N_2 \}_i$. Furthermore, $T_3$ is pronilpotent, and both $T_1$ and $T_2$ normalize $T_3$. Hence, $T_3T_2T_1$ admits a series satisfying the required conditions. Arguing by induction, let $N_j = \gamma_{\infty}(N_{j-1})$ and let $T_j$ be the basis normalizer of the Sylow basis $\{ P_i \cap N_{j-1} \}_i$. We have $G = N_jT_j \cdots T_2T_1$, and the subgroup $T_j \cdots T_2T_1$ has the required series.

 Let $h$ be the Fitting height of $G$.  Then $G = T_h \cdots T_1$ has a series with the required properties.
\end{proof}
\medskip

First, we prove Theorem \ref{MT} (1) in the prosoluble case.

\begin{lemma} \label{solcase}
    Let $G$ be a prosoluble group in which $\pi(G)$ is finite and $|x^G|$ is finite for every $x \in X(G)$. Then the second commutator subgroup $G''$ is finite. 
\end{lemma}

\begin{proof}
        By Proposition \ref{Gvirtnilp}, $G$ is virtually nilpotent. Hence, $G$ has finite Fitting height. Since $\pi(G)$ is finite, Lemma \ref{auxlemma} provides a series
\[
1 \le P_n \le P_nP_{n-1} \le \dots \le P_nP_{n-1} \cdots P_1 = G
\]
    where each $P_i$ is a pro-$p$ subgroup of $G$ for some $p \in \pi(G)$, and $P_i$ normalizes $P_j$ whenever $i \le j$.

    Let $x, y \in G$. Write $x = x_n \cdots x_1$ and $y = y_n \cdots y_1$ where $x_i, y_i \in P_i$ for $i = 1, \dots, n$. Hence,
    \[
    [x, y] = \prod [x_i,y_j]^{\alpha_{ij}}
    \]
    for suitable $\alpha_{ij} \in G$. Note that  $[x_i, y_j] \in P_i$ if $i \ge j$,  whereas  $[x_i, y_j] \in P_j$  if $i \le j$. 
    Hence, $[x_i, y_j] \in X(G) \subseteq \Delta(G)$ for every $i,j = 1, \dots, n $. Since this product is finite and $\Delta(G)$ is a normal (abstract) subgroup, $[x, y] \in \Delta(G)$. Thus, every commutator of $G$ belongs to $\Delta(G)$. From Theorem \ref{known} (1), we conclude that $G''$ is finite. 
\end{proof}
\medskip

For an arbitrary group $G$, not necessarily prosoluble, we will need the following lemma.

\begin{lemma} \label{[H,x]} Let $G$ be a profinite group and let $H$ be the topological closure of $\Delta(G)$. 
Then, for every $x \in \Delta(G)$, the subgroup $[H,x]$ is finite.
\end{lemma}

\begin{proof}
    Throughout the proof, $\Delta = \Delta(G)$ and its subgroups are abstract groups. Since $\Delta$ is an $FC$-group, its commutator subgroup $\Delta '$ is locally finite (see \cite[Section 14.5]{robinson}). Moreover, as $|x^G|$ is finite, $[\Delta, x]$ is finitely generated. Thus $[\Delta, x]$ is finite and hence closed. 
    Since $[H,x]  = [\overline{\Delta}, x] \le \overline{[\Delta, x]}= [\Delta, x],$ we conclude that $[H, x]$ is finite, as claimed. 
\end{proof}
\medskip

Now we are ready to prove Theorem \ref{MT}, which we restate here for the reader's convenience.

\begin{thmMT}  Let $G$ be a profinite group in which $|x^G|$ is finite whenever $x$ is a commutator of prime power order. Assume that $\pi(G)$ is finite. Then
    \begin{enumerate}
\item The second commutator subgroup $G''$ is finite.
\item $G$ has an open normal subgroup of nilpotency class at most $3$.
    \end{enumerate}

\end{thmMT}

\begin{proof}
    By Proposition \ref{Gvirtnilp}, $G$ is virtually nilpotent.
\medskip

    Part (1).
    Since $G$ is virtually nilpotent, its prosoluble radical $R = R(G)$ has finite index in $G$ and is soluble. We argue by induction on $[G:R]$. 
    If $[G:R]=1$, then $G$ is prosoluble and the result follows from Lemma \ref{solcase}.  So we assume that $[G:R] \geq 2$. 

    Let $H=G^{(\infty)}$ be the intersection of all the terms of the derived series of $G$.  Since $G/H$ is prosoluble, Lemma \ref{solcase} implies that $G''/H$ is finite. Therefore, it suffices to prove that $H$ is finite. If $[H:R(H)] < [G:R]$, then, by induction, $H''$ is finite. Since $H''=H$,  we are done. Thus, we may assume that 
    $[H:R(H)] = [G:R]$. Replacing $G$ by $H$, if necessary, we may therefore assume that $G=G'$.
    
    Suppose that $G$ has a normal subgroup $M$ such that $R < M < G$. Since $[M:R(M)] < [G:R]$, it follows by induction that  $M''$ is finite. 
    Moreover, $M/M'' \leq R(G/M'')$. The prosoluble radical of $G/M''$ has index strictly smaller than $[G:R]$. Hence, by induction, $G''/M''$ is finite. Since  $M''$ is also finite, we conclude  that $G'' = G$ is finite.  Therefore, we may assume that no such subgroup $M$ exists, and hence $G/R$ is simple. 

    Since $G=G'$, it follows that $ \langle X(G) \rangle  = G$.
    Taking into account that $G/R$ is a nonabelian simple group, there exists $x\in X(G)$ such that $[G,x]R=G$. It follows that the quotient $G/[G,x]$ is soluble. As $G=G'$, we deduce that $G=[G,x]$. By Lemma \ref{[H,x]},
    $[G,x]$ is finite. Therefore, $G$ is finite,
    as we wanted to prove.
\medskip

    Part (2).
Let $F = F(G)$ be the Fitting subgroup of $G$. 
By Lemma \ref{pronilp}, every commutator of $F$ belongs to $\Delta(G)$. It follows from Theorem \ref{known} (2) that $F$ has an open normal subgroup $L$ of nilpotency class at most $3$. Since $G$ is virtually nilpotent, $[G:L]$ is finite. Let $H$ be the normal core of $L$.
 Thus $[G:H]$ is finite as well. Therefore, $H$ is an open normal subgroup of $G$ of nilpotency class at most $3$, as claimed. 
\end{proof}


\section*{Acknowledgments}  This work was carried out during a visit of the second author to the Department of Mathematics of the University of Padova. The second author thanks the department for its hospitality. P. Shumyatsky was supported by CNPq and FAPDF. E. Detomi is a member of GNSAGA (INDAM).


\bibliographystyle{abbrv}

\end{document}